\documentclass[12pt]{amsart}
\usepackage{amssymb,latexsym,amsmath,amscd,amsthm,graphicx}
\usepackage[all]{xy}
\theoremstyle{plain}

\newtheorem{theorem}[subsection]{Theorem}

\newtheorem{lemma}[subsection]{Lemma}

\theoremstyle{definition}
\newtheorem{prop}[subsection]{Proposition}
\newtheorem{cor}[subsection]{Corollary}

\newtheorem{remark}[subsection]{Remark}

\newtheorem{note}[subsection]{Note}

\newcommand{\uu}{\cup}
\newcommand{\ii}{\cap}
\newcommand{\UU}{\bigcup}
\newcommand{\II}{\bigcap}

\newcommand{\ci}{\subseteq}
\newcommand{\sci}{\subset}
\newcommand{\es}{\emptyset}
\newcommand{\set}[1]{\{#1\}}

\newcommand{\ga}{\alpha}
\newcommand{\gb}{\beta}

\renewcommand{\gg}{\gamma}
\newcommand{\gh}{\eta}

\newcommand{\gk}{\kappa}

\newcommand{\gn}{\nu}
\newcommand{\go}{\omega}

\newcommand{\gs}{\sigma}
\newcommand{\gt}{\tau}
\newcommand{\gx}{\xi}

\newcommand{\gG}{\Gamma}

\newcommand{\gO}{\Omega}

\newcommand{\tbf}{\textbf}
\newcommand{\tit}{\textit}

\newcommand{\D}[1]{\mathbb{#1}}

\newcommand{\te}{\text}

\newcommand{\ol}{\overline}
\newcommand{\ul}{\underline}

\newcommand{\nd}{\noindent}

\newcommand{\vp}{\varphi}
\begin{document}

\title{Quantization dimension for Gibbs-like measures on cookie-cutter sets}
To appear, Kyoto Journal of Mathematics, please see the published version.

\author{ Mrinal Kanti Roychowdhury}
\address{Department of Mathematics\\
The University of Texas-Pan American\\
1201 West University Drive\\
Edinburg, TX 78539-2999, USA.}
\email{roychowdhurymk@utpa.edu}


\date{}
\maketitle

\pagestyle{myheadings}\markboth{Mrinal Kanti Roychowdhury}{Quantization dimension for Gibbs-like measures on cookie-cutter sets}
\begin{abstract}
In this paper using Banach limit we have determined a Gibbs-like measure $\mu_h$ supported by a cookie-cutter set $E$ which is generated by a single cookie-cutter mapping $f$.
For such a measure $\mu_h$ and $r\in (0, +\infty)$ we have shown that there exists a unique $\gk_r \in (0, +\infty)$ such that $\gk_r$ is the quantization dimension function of the probability measure $\mu_h$, and established its functional
relationship with the temperature function of the thermodynamic formalism.  The temperature function is commonly used to perform the
multifractal analysis, in our context of the measure $\mu_h$. In addition, we have proved that the $\gk_r$-dimensional lower quantization coefficient of order $r$ of the probability measure is positive.
\end{abstract}
\section{Introduction}

Quantization dimension is one of the most important objects in the quantization problem, which has a deep background in information theory and engineering technology (cf. \cite{BW, GG, GN, Z}). It characterizes in a natural way the
asymptotic property of the error when approximating a given probability measure by a
discrete probability measure of finite support in the sense of $L_r$-metrics. Given a Borel
probability measure $\mu$ on $\D R^d$, a number $r \in (0, +\infty)$
and a natural number $n \in \D N$, the $n$th \tit{quantization
error} of order $r$ of $\mu$, is defined by
\[V_{n, r}(\mu):=\te{inf}\set{\int d(x, \ga)^r d
\mu(x) : \ga \sci \D R^d, \, \te{card}(\ga) \leq n},\]
where $d(x, \ga)$ denotes the distance from the point $x$ to the set
$\ga$ with respect to a given norm $\|\cdot\|$ on $\D R^d$. Note
that if $\int \| x\|^r d\mu(x)<\infty$ then there is some set $\ga$ for
which the infimum is achieved (cf. \cite{GL1}). The set $\ga$ for which the infimum is achieved is called an optimal set of $n$-means or $n$-optimal set of order $r$ for $0<r<+\infty$.
The upper and the lower
quantization dimensions of order $r$ of $\mu$ are defined to be
\[\ol D_r(\mu): =\limsup_{n \to \infty} \frac{r\log n}{-\log V_{n, r}(\mu)}; \ \ul D_r(\mu): =\liminf_{n \to \infty} \frac{r\log n}{-\log V_{n, r}(\mu)}. \]
If $\ol D_r(\mu)$ and $\ul D_r(\mu)$ coincide, we call the common value the \tit{quantization dimension of order $r$} of the probability measure $\mu$, and is
denoted by $D_r:=D_r(\mu)$. For $s>0$, we define the $s$-dimensional upper and lower quantization coefficients of order $r$ of $\mu$ by $\limsup_n n V_{n,r}^{\frac s r} (\mu)$ and $\liminf_n n V_{n,r}^{\frac s r} (\mu)$ respectively. One sees that the quantization dimension is actually a function $r\mapsto D_r$ which measures the asymptotic rate at which $V_{n, r}$ goes to zero. If $D_r$ exists, then one can write
\[\log V_{n, r} \sim \log [(\frac{1}{n})^{r/D_r}].\]
For probabilities with non-vanishing absolutely continuous part the numbers $D_r$ are all equal to the dimension $d$ of the underlying space, but for singular probabilities the family $(D_r)_{r>0}$ gives an interesting description of their geometric (multifractal) structures.

Let $S_1, S_2, \cdots, S_N$ be contractive similitudes from $\D R^d$ into itself, where $N\geq 2$ is a positive integer. Let $s_i$ be the contraction ratio of $S_i$ for all $1\leq i\leq N$. Then for a given probability vector $(p_1, p_2, \cdots, p_N)$ there exists a unique Borel probability measure $\mu$ (cf. \cite{H}) satisfying the condition
\[\mu=\sum_{j=1}^N p_j \mu \circ S_j^{-1}.\]
Let the iterated function system $\set{S_1, S_2, \cdots, S_N}$ satisfy the open set condition: there exists a bounded nonempty open set $U \sci \D R^d$ such that $\UU_{j=1}^N S_j(U) \sci U$ and $S_i(U) \II S_j(U) =\es$ for $1\leq i\neq j\leq N$. The iterated function system satisfies the strong open set condition if $U$ can be chosen such that $U\ii J \neq \es$, where $J$ is the limit set of the iterated function system. Under the open set condition,  Graf and Luschgy showed that the quantization dimension function $D_r:=D_r(\mu)$ of the probability measure $\mu$ exists, and satisfies the following relation (cf. \cite{GL1, GL2}):
\begin{equation*} \sum_{j=1}^N(p_j s_j^r)^{{\frac{D_r}{r+D_r}}}=1.\end{equation*}
In fact, they proved a stronger result namely that the quantization dimension $D_r$ also satisfies (cf. \cite{GL3}):
\begin{equation} \label{eq1} 0<\liminf_n n V_{n,r}^{\frac{D_r}{r}}(\mu)\leq \limsup_n n V_{n,r}^{\frac{D_r}{r}}(\mu)<+\infty.\end{equation}
Under the open set condition, Lindsay and
Mauldin (cf. \cite{LM}) determined the quantization dimension for an $F$-conformal measure $m$ associated
with a conformal iterated function system determined by finitely
many conformal mappings. They established a relationship between the
quantization dimension and the temperature function of the thermodynamic formalism arising in multifractal analysis. Later, in \cite{R1} the author studied the quantization dimension of Moran measures on the Moran sets of which potential functions are defined in terms of the similarity ratios
and probability vectors. Then the quantization problem was solved in \cite{R2} for the image
measures of Bowen's Gibbs measures supported by the one-sided Bernoulli shifts under
the coding maps on the Moran sets. Mixed quantization dimension and its relationship with the temperature function was studied in \cite{WD} by Wang and Dai. But from the work in \cite{LM, R1, R2, WD} it was not known whether the $D_r$-dimensional lower quantization coefficient is positive, i.e., whether $\liminf_n n V_{n,r}^{\frac {D_r} r} (\mu)>0$, where $D_r:=D_r(\mu)$ is the quantization dimension of the probability measure $\mu$.

In this paper, using Banach limit we have defined a Gibbs-like measure $\mu_h$ supported by a cookie-cutter set $E$, where $E$ is the limit set generated by a cookie-cutter mapping $f$
and $h:=\te{dim}_{\te{H}}(E)$ is the Hausdorff dimension of the set $E$ (cf. \cite{F1}).  For this measure $\mu_h$ we have shown that for each $r\in (0, +\infty)$ there
exists a unique $\gk_r \in (0, +\infty)$ such that
\begin{equation} \label{eq2000} \lim_{k\to \infty}\frac 1 k\log\sum_{\gs \in \gO_k} \left(\mu_h(J_\gs) \|\vp_\gs'\|^r\right)^{\frac{\gk_r}{r+\gk_r}}=0,  \end{equation}
and the above $\gk_r$ is the quantization dimension $D_r:=D_r(\mu_h)$ of order $r$ of the probability measure $\mu_h$.
It is known that, the
singularity exponent $\gb(q)$ (also known as the temperature function)
satisfies the usual equation
\begin{equation}\label{eq2001}  \lim_{k\to\infty}\frac 1 k\log
  \sum_{\gs \in \gO_k}\left(\mu_h(J_\gs)\right)^q \|\vp_\gs'\|^{\gb(q)}=0, \end{equation}
and that the spectrum $f(\ga)$ is the Legendre transform of
$\gb(q)$. Comparing \eqref{eq2000} and \eqref{eq2001}, we see that if
$q_r=\frac{D_r}{r+D_r}$, then $\gb(q_r)=rq_r$, that is, the
quantization dimension function of order $r$ of the probability measure $\mu_h$
has a relationship with the temperature function of the thermodynamic
formalism arising in multifractal analysis (for thermodynamic
formalism, multifractal analysis and the Legendre transform one could
see \cite{F1}). The significant difference of the work in this paper and the work in \cite{LM, R1, R2, WD} is that, in addition to determine the quantization dimension function $\gk_r$ and its relationship with the temperature function of a probability measure, by Proposition~\ref{prop6} in this paper, we have proved a stronger result namely that the $\gk_r$-dimensional lower quantization coefficient $\liminf_n n V_{n,r}^{\frac {\gk_r} r} (\mu_h)$ of order $r$ of the probability measure $\mu_h$ is positive. Quantization problem for a general probability measure is still open.

\section{Basic definitions, lemmas and propositions}

In this paper, $\D R^d$ denotes the $d$-dimensional Euclidean space equipped with a metric $d$ compatible with the Euclidean topology.
Let us write
\begin{align*}  V_{n, r}(\mu):&=\te{inf}\set{\int d(x, \ga)^r d\mu(x) : \ga \sci \D R^d, \te{ card}(\ga)\leq n},\\
u_{n, r}(\mu):&=\te{inf}\set{\int d(x, \ga\uu U^c)^r d\mu(x) : \ga \sci \D R^d, \te{ card}(\ga)\leq n},
 \end{align*} where $U$ is a set which comes from the open set condition and $U^c$ denotes the complement of $U$. We see that
\[u_{n, r}^{1/r} \leq V_{n, r}^{1/r}:=e_{n, r}.\]
We call sets $\ga_n \sci \D R^d$, for which the above infimums are achieved, $n$-optimal sets for $e_{n, r}, V_{n, r}$ or $u_{n, r}$ respectively. As stated above, Graf and Luschgy have shown that $n$-optimal sets exist when $\int \|x\|^r d\mu(x)<\infty$.

\subsection{Cookie-cutter set:} A mapping $f$ is called a cookie-cutter, if there  exists a finite collection of disjoint closed intervals $J_1, J_2, \cdots, J_N \sci J =[0, 1]$, such that

$(C1)$ $f$ is defined in a neighborhood of each $J_j$, $1\leq j\leq N$, the restriction of $f$ to each initial interval $J_j$ maps $J_j$ bijectively onto $J$, and the corresponding
branch inverse is denoted by $\vp_j:=(f|_{J_j})^{-1} : J \to J_j$;

$(C2)$ $f$ is differentiable with H\"older continuous derivative $f'$, i.e., there exist constants $c>0$ and $\gg \in (0, 1]$ such that for $x, y\in J_j$,
$1\leq j\leq N$,
\[\left|f'(x)-f'(y)\right|\leq c|x-y|^{\gg};\]

$(C3)$ $f$ is boundedly expanding in the sense that there exist constants $b$ and $B$
\[1<b:=\inf_x\set{|f'(x)|} \leq \sup_x\set{|f'(x)| }:=B<+\infty.\]
$[\UU_{j=1}^N J_j; c, \gg, b, B]$ is called the defining data of the cookie-cutter mapping $f$.
Let $\gO_0$ be the empty set. For $n\geq 1$, define
\[\gO_{n}=\set{1, 2, \cdots, N}^n, \ \gO_\infty=\set{1, 2, \cdots, N}^{\D N} \te{ and } \gO=\UU_{k=0}^\infty \gO_k.\]
Elements of $\gO$ are called words. For any $\gs \in \gO$ if $\gs=(\gs_1, \gs_2, \cdots, \gs_n)  \in \gO_n$, we write $\gs^-=(\gs_1, \gs_2, \cdots, \gs_{n-1})$ to denote the word obtained by deleting the last letter of $\gs$,  $|\gs|=n$ to denote the length of $\gs$, and $\gs|_k:=(\gs_1, \gs_2, \cdots, \gs_k)$, $k\leq n$, to denote the truncation of $\gs$ to the length $k$.  For any two words  $\gs=(\gs_1, \gs_2, \cdots, \gs_k)$ and $\gt=(\gt_1, \gt_2, \cdots, \gt_m)$, we write $\gs\gt=\gs\ast \gt=(\gs_1,  \cdots, \gs_k, \gt_1,  \cdots, \gt_m)$ to denote the juxtaposition of $\gs, \gt \in \gO$.  A word of length zero is called the empty word and is denoted by $\es$. For $\gs \in \gO$ and $\gt\in \gO\uu \gO_\infty$ we say $\gt$ is an extension of $\gs$, written as $\gs\prec \gt$, if $\gt|_{|\gs|}=\gs$. For $\gs \in \gO_k$, the cylinder set $C(\gs)$ is defined as $C(\gs)=\set{\gt  \in \gO_\infty : \gt|_k =\gs}$.
For $\gs=(\gs_1, \gs_2, \cdots, \gs_n) \in \gO_n$, let us write $\vp_\gs=\vp_{\gs_1} \circ \cdots \circ\vp_{\gs_n}$, and define the rank-$n$ basic interval  corresponding to $\gs$ by
\[J_\gs=J_{(\gs_1, \gs_2, \cdots, \gs_n)}=\vp_{\gs}(J),\]
where $1\leq \gs_k \leq N$, $1\leq k\leq n$. If $\gs=\es$, then we identify $\vp_\es$ to be the identity mapping on $J$ and write $J_\es=J$. By $|J_\gs|$ we mean the
diameter of the set $J_\gs$ for all $\gs \in \gO$. It is easy to see that the set of basic intervals $\set{J_\gs : \gs \in \gO}$ has the following net properties:

$(i)$ $J_{\gs \ast j} \sci J_\gs$ for each $\gs \in \gO_n$ and $1\leq j\leq N$ for all $n\geq 1$;

$(ii)$ $J_{\gs} \II J_{\gt}=\es$, if $\gs, \gt \in \gO_n$ for all $n\geq 1$ and $\gs \neq \gt$.

Since $\vp_{j}$ is a branch inverse of $f$, where $1\leq j\leq N$, for all $x \in J$, we have $f(\vp_{j}(x))=x$, and so $|f'(\vp_j(x))|\cdot|\vp_j'(x)|=1$, which yields
\begin{equation} \label{eq1} B^{-1} \leq \left|\vp_{j}'(x)\right|\leq b^{-1}. \end{equation}
Choose $x, y$ to be the end points of $J$, and then $\vp_\gs(x), \vp_\gs(y)$ are the end points of $J_\gs$ for each $\gs \in \gO$, and so by mean value theorem, we have  \[|J_\gs|=|\vp_\gs(x)-\vp_\gs(y)|=|\vp_\gs'(w)||x-y|=|\vp_\gs'(w)|\]
 for some $w\in J_\gs$. Thus $B^{-n} \leq |J_\gs| \leq b^{-n}$ for any $\gs \in \gO_n$, and thus the diameter $\left|J_\gs\right| \to 0$ as $|\gs|\to \infty$. Since given $\gs =(\gs_i)_{i=1}^\infty \in \gO_\infty$ the diameters of the compact sets $J_{\gs|_k}, \, k\geq 1,$ converge to zero and since they form a descending family, the set
 \[\II_{k=0}^\infty J_{\gs|_k}\]
 is a singleton and therefore, if we denote its element by $\pi(\gs)$, this defines the coding map $\pi : \gO_\infty \to J$. The main object of our interest is the limit set
 \[E:=\pi(\gO_\infty)=\UU_{\gs \in \gO_\infty} \II_{k=0}^\infty J_{\gs|_k}.\]
Moreover, $\pi(C(\gs))=E\ii J_\gs$ for $\gs\in \gO$. With the net properties it follows that $E$ is a perfect, nowhere dense and totally disconnected subset of $J$.
The set $E$ is called the \tit{cookie-cutter} set.

Let $\ell^\infty$ be the set of all bounded sequences $x=(x_n)_{n\in \D N}$ of real or complex numbers, which form a vector space with respect to point-wise addition and multiplication by a scalar. It is equipped with the norm $\|x\|=\sup_n{|x_n|}$. The normed space $\ell^\infty$ is complete with respect to the metric $\|x-y\|$, and so it forms a Banach space. By the Hahn-Banach theorem (cf. \cite[p. 102-104]{Y}), there exists a linear functional $L : \ell^\infty \to \D R$ for which

$(i)$ $L$ is linear;

$(ii)$ $L((x_n)_{n\in \D N})=L((x_{n+1})_{n\in \D N})$;

$(iii)$ $\liminf_{n\to \infty} (x_n)\leq L((x_n)_{n\in\D N})\leq \limsup_{n\to \infty} (x_n)$.

\nd The functional $L$, defined above, is called a Banach limit.  The use of the Banach limit is rather a standard tool in producing an invariant measure from a given measure.

Let us now prove the following lemma.
\begin{lemma}  \label{lemma1111} (\tbf{Bounded variation principle})
There exists a constant $1<\gx<+\infty$ such that for each $\gs \in \gO_{n}$, and $x, y \in J_\gs$, we have
\[\gx^{-1} \leq \frac{|(f^n)'(x)|}{|(f^n)'(y)|} \leq \gx,\]
where $f^n= f \circ f \circ\cdots\circ f$ represents the n-fold composition of $f$ with itself.
\end{lemma}
\begin{proof}
Note that for each $k\leq n$ and $\gs=(\gs_1, \gs_2, \cdots, \gs_n)\in \gO_{n}$,  $f^{k-1}$ maps $J_\gs$ diffeomorphically to the set $\vp_{\gs_k} \circ \vp_{ \gs_{k+1}}\circ \cdots \circ \vp_{\gs_n}(J)$, and so
\[\left|f^{k-1}(x)-f^{k-1}(y)\right| \leq \te{diam}\left( \vp_{\gs_k} \circ \vp_{\gs_{k+1}}\circ \cdots \circ \vp_{\gs_n}(J)\right)=|\vp_{\gs_k} \circ \vp_{ \gs_{k+1}}\circ \cdots \circ \vp_{\gs_n}(J)|.\]
By mean value theorem,
\begin{align*}
&|\vp_{\gs_{k}} \circ \vp_{\gs_{k+1}}\circ \cdots \circ \vp_{\gs_n}(J)|\\
&=\sup_{x, y\in J} \left|\vp_{\gs_k} \left(\vp_{\gs_{k+1}}\circ \cdots \circ \vp_{\gs_n}(x)\right)-\vp_{\gs_k} \left(\vp_{\gs_{k+1}}\circ \cdots \circ \vp_{ \gs_n}(y)\right)\right|\\
&\leq b^{-1} \left|\vp_{\gs_{k+1}}\circ \cdots \circ \vp_{\gs_n}(J)\right|.
\end{align*}
Thus proceeding inductively,
\[\left |\vp_{\gs_k} \circ \vp_{\gs_{k+1}}\circ \cdots \circ \vp_{\gs_n}(J)\right |\leq b^{-(n-k+1)}.\]
Then, H\"older continuity of $f'$ gives
\begin{align*}\left  |f'(f^{k-1}(x))-f'(f^{k-1}(y))\right | \leq c \left |f^{k-1}(x)-f^{k-1}(y)\right|^\gg\leq cb^{-(n-k+1)\gg},\end{align*}
and so by mean value theorem and the assumption $|f'|>1$, we have
\begin{align*} &\Big |\log |f'(f^{k-1}(x))|-\log |f'(f^{k-1}(y))| \Big |\\
&\leq \Big ||f'(f^{k-1}(x))| -|f'(f^{k-1}(y))| \Big |\\
& \leq cb^{-(n-k+1)\gg}.
\end{align*}
Therefore, by the above inequality and the chain rule,
\begin{align*}
&\Big |\log |(f^n)'(x)|-\log |(f^n)'(y)|\Big |\\
&=\Big |\sum_{k=1}^n \log |f'(f^{k-1}(x))|-\sum_{k=1}^n \log |f'(f^{k-1}(y))|\Big |\\
&\leq \sum_{k=1}^n \Big |\log |f'(f^{k-1}(x))|-\log |f'(f^{k-1}(y))|\Big |\\
&\leq  \sum_{k=1}^n cb^{-(n-k+1)\gg}\leq \frac{cb^{-\gg}}{1-b^{-\gg}}.
\end{align*}
Take $\gx=\exp\left\{\frac{c}{b^{\gg}-1}\right\}$. Since $\frac{c}{b^{\gg}-1}>0$, we have $1< \gx<+\infty$, and thus the lemma follows.
\end{proof}

Let us now prove the following proposition.

\begin{prop}. \label{prop31} (\tbf{Bounded distortion principle})  For any $n\geq 1$, $\gs\in \gO_n$, $x \in J_\gs$, we have
\[\gx^{-1} \leq  |(f^n)'(x)|\cdot |J_\gs| \leq \gx.\]
Moreover, for each $1\leq j\leq N$, we get $\left|J_{\gs\ast j}\right| \geq \gx^{-2} B^{-1}|J_\gs|$, where $\gx$ is the constant of Lemma~\ref{lemma1111}.
\end{prop}

\begin{proof}
Note that for $\gs \in \gO_n$,  $f^n: J_\gs \to J$ is a differentiable bijection. So by mean value theorem, if $y, z \in J_\gs$, there exists $w \in J_\gs$ such that
\[f^n(y)-f^n(z)=(f^n)'(w) (y-z).\]
Choose $y, z$ to be the end points of $J_\gs$, and then $f^n(y), f^n(z)$ are the end points of $J$, and so
\begin{equation*} \label{eq1111}|J|=|(f^n)'(w)|\cdot|J_\gs|, \te{ i.e., } |(f^n)'(w)|\cdot|J_\gs|=1.\end{equation*}
Hence, using bounded variation principle, we have
\begin{equation} \label{eq1112} \gx^{-1} \leq  |(f^n)'(x)|\cdot |J_\gs| \leq \gx\end{equation}  for all $x \in J_\gs$.
Now let $1\leq j\leq N$ and  $x \in J_{\gs\ast j}$. Then using \eqref{eq1112}, we have
\[\gx^{-1} \leq  |(f^{n+1})'(x)|\cdot |J_{\gs\ast j} |=|(f'(f^n(x))|\cdot |(f^n)'(x)|\cdot |J_{\gs\ast j}| \leq B |(f^n)'(x)| \cdot |J_{\gs\ast j}|.\]  Since $J_{\gs\ast j} \ci J_\gs$, we have $x \in J_\gs$. Hence using \eqref{eq1112} again, we have
\[ |J_{\gs\ast j}|\geq \gx^{-2} B^{-1}|J_\gs|.\]
Thus the proof of the proposition is yielded.
\end{proof}

 \begin{prop} \label{prop32}
For any $n\geq 1$, let $\gs\in\gO_{n}$, and  $x, y \in J$. Let $\gx$ be the constant of Lemma~\ref{lemma1111}. Then,
\[\gx^{-1} |\vp_\gs'(y)| \leq |\vp_\gs'(x)| \leq \gx |\vp_\gs'(y)|.\]
\end{prop}
\begin{proof}
For $\gs\in \gO_{n}$ and $x \in J$, we know $f^n(\vp_\gs(x))=x$. Thus
\[|(f^n)'(\vp_\gs(x))|\cdot|\vp_\gs'(x)|=1.\]
Again for all $x \in J$, $\vp_\gs(x) \in J_\gs$. Hence, Lemma~\ref{lemma1111} yields
\[\gx^{-1} |\vp_\gs'(y)| \leq |\vp_\gs'(x)| \leq \gx |\vp_\gs'(y)|,\]
and thus the proposition is obtained.
\end{proof}

Let us now prove the following two lemmas.
\begin{lemma} \label{lemma345}
Let $\gs, \gt\in \gO$. Then
\[\gx^{-1} \|\vp_\gs'\|\|\vp_{\gt}'\|\leq \|\vp_{\gs\gt}'\|\leq \|\vp_\gs'\|\|\vp_{\gt}'\|, \]
$\gx$ is the constant of Lemma~\ref{lemma1111}
\end{lemma}
\begin{proof}
For any $x\in J$, we have $|\vp_{\gs\gt}'(x)|=|\vp_\gs'(\vp_\gt(x))|\cdot|\vp_\gt'(x)|$, and thus by Proposition~\ref{prop32}, for any $y \in J$,
\[\gx^{-1}|\vp_\gs'(y)|\cdot |\vp_\gt'(x)| \leq |\vp_\gs'(\vp_\gt(x))|\cdot|\vp_\gt'(x)|=|\vp_{\gs\gt}'(x)|\leq \|\vp_\gs'\|\cdot \|\vp_\gt'\|,\]
and thus the lemma follows.
\end{proof}

\begin{lemma} \label{lemma32}
Let $\gs \in \gO$ and $x \in J$. Then
\[\gx^{-1} |J_\gs| \leq |\vp_\gs'(x)| \leq \gx |J_\gs|,\]
where $\gx$ is the constant of Lemma~\ref{lemma1111}.
\end{lemma}
\begin{proof} Let $x \in J$, and then $\vp_\gs(x) \in J_\gs$ for $\gs \in \gO_n$, $n\geq 1$. We know $f^n(\vp_\gs(x))=x$, and so $|(f^n)'(\vp_\gs(x))|\cdot |\vp_\gs'(x)|=1$. Now use  Proposition~\ref{prop31} to obtain the lemma. \end{proof}

Let us now prove the following lemma.
\begin{lemma}\label{lemma33}
Let $\gs, \gt \in \gO$. Then
\[\gx^{-3} |J_\gs||J_\gt| \leq |J_{\gs\gt}| \leq \gx^3 |J_\gs||J_\gt|,\]
where $\gx$ is the constant of Lemma~\ref{lemma1111}.
\end{lemma}
\begin{proof}
For $\gs, \gt \in \gO$, we have  $|\vp_{\gs\gt}'(x)|=|\vp_\gs'(y)||\vp_\gt'(x)|$ where $y=\vp_\gt(x)$, and $x \in J$. Again by Proposition~\ref{prop32}, for any $x, y \in J$, we have  \[\gx^{-1} |\vp_\gs'(y)| \leq |\vp_\gs'(x)| \leq \gx |\vp_\gs'(y)|.\]
Hence, Lemma~\ref{lemma32} implies
\[\gx^{-3} |J_\gs||J_\gt|\leq \gx^{-1} |\vp_\gs'(y)||\vp_\gt'(x)|=\gx^{-1} |\vp_{\gs\gt}'(x)|\leq |J_{\gs\gt}|\leq \gx |\vp_{\gs\gt}'(x)|\leq \gx^{3} |J_\gs||J_\gt|,\]
and thus the lemma is obtained.
\end{proof}
By Lemma~\ref{lemma345} and the standard theory of sub-additive sequences, the function $Q(t)$ given by
\[Q(t)=\lim_{k\to \infty} \frac 1 k\log  \sum_{\gs\in \gO_k} \|\vp_\gs'\|^t,\]
for any real $t$ exists. It is easy to see that the function $Q(t)$ is strictly decreasing convex and hence continuous in $t$.

\begin{lemma} \label{lemma12}
There exists a unique $h \in (0, +\infty)$ such that $Q(h)=0$.
\end{lemma}

\begin{proof}
Since the function $Q(t)$ is strictly decreasing and continuous on $\D R$, there exists a unique $h \in \D R$ such that $Q(h)=0$. Note that
\[
Q(0)=\lim_{k\to\infty} \frac 1 k \log \sum_{\gs\in \gO_k} 1\geq \lim_{k\to\infty} \frac 1 k \log N^k=\log N \geq \log 2>0.
\]
In order to conclude the proof it therefore suffices to show that $\lim_{t \to +\infty} Q(t)=-\infty$. For $t>0$,
\[Q(t)=\lim_{k\to\infty} \frac 1 k \log \sum_{\gs\in \gO_k}\|\vp_\gs'\|^t \leq \lim_{k\to\infty} \frac 1 k\log  \sum_{\gs\in \gO_k} b^{-kt}=\lim_{k\to\infty} \frac 1 k\log N^{k} -t \log b=\log N-t\log b.\]
Since $b>1$, it follows that $\lim_{t\to +\infty}Q(t)=-\infty$, and hence the lemma follows.
\end{proof}

\begin{note} Lemma~\ref{lemma32} implies $\gx^{-1} |J_\gs| \leq \sup_{x \in J}|\vp_\gs'(x)|=\|\vp_\gs'\| \leq \gx |J_\gs|,$ and so the topological pressure $Q(t)$ can be written as follows:
\[Q(t)=\lim_{k\to \infty} \frac 1 k \log \sum_{\gs\in \gO_k}|J_\gs|^t. \]
The unique $h \in (0, +\infty)$  given by Lemma~\ref{lemma12} is the Hausdorff dimension $\te{dim}_{\te{H}}(E)$ of the cookie-cutter set $E$ (cf. \cite{F1}).

\end{note}

Let us now prove the following proposition, which plays a vital role in the paper.
\begin{prop} \label{prop21}
Let $h \in (0, +\infty)$ be unique such that $Q(h)=0$, and let $s_\ast$ and $s^\ast$ be any two arbitrary real numbers with $0<s_\ast<h<s^\ast$. Then for all $n\geq 1$,
\[\gx^{-3s_\ast}< \sum_{\gs \in \gO_n}|J_\gs|^{s_\ast} \te{ and } \sum_{\gs \in \gO_n}|J_\gs|^{s^\ast}<\gx^{3s^\ast},\]
where $\gx$ is the constant of Lemma~\ref{lemma1111}.
\end{prop}
\begin{proof}
Let $s_\ast<h$. As the pressure function $Q(t)$ is strictly decreasing, $Q(s_\ast)>Q(h)=0$. Then for any positive integer $n$, by Lemma~\ref{lemma33}, we have
\begin{align*}
0&<Q(s_\ast) =\lim_{p\to \infty} \frac 1{np} \log \sum_{\go\in \gO_{np}}|J_\go|^{s_\ast}\leq \lim_{p\to \infty} \frac 1{np} \log \gx^{3(p-1)s_\ast} \left(\sum_{\gs\in \gO_{n}}|J_\gs|^{s_\ast}\right)^p,
\end{align*}
which implies \[0< \frac 1n \log \left(\gx^{3s_\ast} \sum_{\gs\in \gO_{n}}|J_\gs|^{s_\ast}\right) \te{ and so } \sum_{\gs\in \gO_{n}}|J_\gs|^{s_\ast}>\gx^{-3s_\ast}. \]
Now if $h<s^\ast$, then $Q(s^\ast)<0$ as $Q(t)$ is strictly decreasing.  Then for any positive integer $n$, by Lemma~\ref{lemma33}, we have
\begin{align*}
0>Q(s^\ast) =\lim_{p\to \infty} \frac 1{np} \log \sum_{\go\in \gO_{np}}|J_\go|^{s^\ast}\geq \lim_{p\to \infty} \frac 1{np} \log \gx^{-3(p-1)s^\ast} \left(\sum_{\gs\in \gO_{n}}|J_\gs|^{s^\ast}\right)^p,
\end{align*}
which implies \[0>\frac 1n \log \left(\gx^{-3s^\ast}\sum_{\gs\in \gO_{n}}|J_\gs|^{s^\ast}\right) \te{ and so } \sum_{\gs\in \gO_{n}}|J_\gs|^{s^\ast}<\gx^{3s^\ast}. \]
Thus the proposition is obtained.
\end{proof}
\begin{cor} \label{cor1}
Since  $s_\ast$ and $s^\ast$ be any two arbitrary real numbers with $0<s_\ast<h<s^\ast$, from the above proposition it follows that for all $n\geq 1$,
\[\gx^{-3h}\leq \sum_{\gs \in \gO_n}|J_\gs|^h\leq \gx^{3h}.\]

\end{cor}

Let us now prove the following proposition.

\begin{prop} \label{prop11} \tbf{(Gibbs-like measure)}
Let $h\in (0, +\infty)$ be such that $Q(h)=0$. Then there exists a constant $\gh>1$ and a probability measure $\mu_h$ supported by $E$ such that for any $\gs\in \gO$,
\[\gh^{-1} |J_{\gs}|^h\leq \mu_h(J_{\gs}) \leq \gh |J_{\gs}|^h.\]
\end{prop}
\begin{proof}
For $\gs\in \gO$, $n\geq 1$, define
\[\gn_n(C(\gs))=\frac{\sum_{\gt\in \gO_n} (\te{diam} J_{\gs\gt})^h}{\sum_{\gt\in \gO_{|\gs|+n}} (\te{diam}J_{\gt})^h}.\]
Then using Lemma~\ref{lemma33} and Corollary~\ref{cor1}, we have
\[\gn_n(C(\gs))\leq \frac{\gx^{3h}  (\te{diam}J_\gs)^h \sum_{\gt\in \gO_n}(\te{diam} J_\gt)^h}{\sum_{\gt\in \gO_{|\gs|+n}} (\te{diam}J_{\gt})^h}\leq \gx^{9h}(\te{diam}J_\gs)^h, \]
and similarly, $\gn_n(C(\gs))\geq \gx^{-9h}(\te{diam}J_\gs)^h$. Thus for a given $\gs \in \gO$, $\set{\gn_n(C(\gs))}_{n=1}^\infty$ is a bounded sequence of real numbers,
and so Banach limit, denoted by Lim, is defined. For $\gs \in \gO$, let
\[\gn(C(\gs))=\te{Lim}_{n\to \infty} \gn_n(C(\gs)).\]
Then
\[\sum_{j=1}^N \gn(C(\gs  j))=\te{Lim}_{n\to \infty} \sum_{j=1}^N \frac{\sum_{\gt\in \gO_n} (\te{diam} J_{\gs j \gt})^h}{\sum_{\gt\in \gO_{|\gs|+1+n}} (\te{diam}J_{\gt})^h}=\te{Lim}_{n\to \infty} \frac{\sum_{\gt\in \gO_{n+1}} (\te{diam} J_{\gs \gt})^h}{\sum_{\gt\in \gO_{|\gs|+n+1}} (\te{diam}J_{\gt})^h}, \]
and so
\[\sum_{j=1}^N \gn(C(\gs  j))=\te{Lim}_{n\to \infty} \gn_{n+1}(C(\gs))=\te{Lim}_{n\to \infty} \gn_{n}(C(\gs))=\gn(C(\gs)).\]
Thus by Kolmogorov's extension theorem, $\gn$ can be extended to a unique Borel probability measure $\gg$ on $\gO_\infty$. Let $\mu_h$ be the image measure of $\gg$
under the coding map $\pi$, i.e., $\mu_h=\gg\circ \pi^{-1}$. Then $\mu_h$ is a unique Borel probability measure supported by $E$. Moreover, for any $\gs\in \gO$,
\[\mu_h(J_\gs)=\gg(C(\gs))=\te{Lim}_{n\to \infty} \gn_n(C(\gs)) \leq \te{Lim}_{n\to\infty}  \gx^{9h}(\te{diam}J_\gs)^h= \gx^{9h}(\te{diam}J_\gs)^h,\]
and similarly,
 \[\mu_h(J_\gs)\geq \gx^{-9h}(\te{diam}J_\gs)^h.\]
Write $\gh=\gx^{9h}$, and then $\gh>1$, and thus the proof of the proposition is complete.
\end{proof}

For the above measure $\mu_h$, known as \tit{Gibbs-like measure}, we will determine the quantization dimension function and its functional relationship with the temperature function of the thermodynamic formalism.

Let us now prove the following lemma.
\begin{lemma} \label{lemma3}    Let $\mu_h$ be the Gibbs-like measure as defined in Proposition~\ref{prop11}. Then there exists a constant $L>1$ such that for $\gs, \gt\in \gO$,
\[L^{-1} \mu_h(J_\gs)\mu_h(J_\gt)\leq \mu_h(J_{\gs\gt})\leq L \mu_h(J_\gs)\mu_h(J_\gt).\]
\end{lemma}
\begin{proof}
Let $\gs, \gt \in \gO$. Then by Lemma~\ref{lemma33} and Proposition~\ref{prop11}, we have
\[\mu_h(J_{\gs\gt}) \leq \gh|J_{\gs\gt}|^h\leq \gh\gx^{3h}|J_\gs|^h|J_\gt|^h\leq \gh^3\gx^{3h} \mu_h(J_\gs)\mu_h(J_\gh),\]
and similarly, $\mu_h(J_{\gs\gt})\geq \gh^{-3}\gx^{-3h} \mu_h(J_\gs)\mu_h(J_\gh)$. Take $L=\gh^3 \gx^{3h}$. As $h>0$, $\gx>1$ and $\gh>1$, it follows that $L>1$, and thus
\[L^{-1}  \mu_h(J_\gs)\mu_h(J_\gh)\leq \mu_h(J_{\gs\gt})\leq L\mu_h(J_\gs)\mu_h(J_\gh),\]
which is the lemma.
\end{proof}

\subsection{Topological pressure:} For $q, t\in \D R$ and $n\geq 1$, let us write
\[Z_n(q, t)=\sum_{\gs\in \gO_n}\left(\mu_h(J_\gs)\right)^q \|\vp_\gs'\|^t.\] Then for $n, p\geq 1$,
\[Z_{n+p}(q, t)=\sum_{\gs\in \gO_n} \sum_{\gt\in \gO_p}\left(\mu_h(J_{\gs\gt})\right)^q \|\vp_{\gs\gt}'\|^t .\] Let us first assume $q\geq 0$. Then by Lemma~\ref{lemma345} and Lemma~\ref{lemma3}, it follows that if $t\geq 0$ then
\[ Z_{n+p}(q, t) \leq L^q Z_n(q, t)Z_p(q, t),\]
and if $t<0$ then
\[ Z_{n+p}(q, t) \leq L^{q}\gx^{-t} Z_n(q, t)Z_p(q, t).\]
Let us now assume $q<0$. Using the same argument, if $t\geq 0$ then
\[ Z_{n+p}(q, t) \leq L^{-q} Z_n(q, t)Z_p(q, t),\]
and if $t<0$ then
\[ Z_{n+p}(q, t) \leq L^{-q}\gx^{-t} Z_n(q, t)Z_p(q, t).\]
Hence by the standard theory of sub-additive sequences, $\lim_{k\to\infty} \frac 1 k \log Z_k(q, t)$ exists (cf. \cite[Corollary 1.2]{F1}). Let us denote it by $P(q, t)$, i.e.,
\begin{equation} \label{eq2} P(q, t)=\lim_{k\to\infty} \frac 1 k \log \sum_{\gs\in \gO_k}\left(\mu_h(J_\gs)\right)^q \|\vp_\gs'\|^t.\end{equation}

The following proposition states the well-known properties of the function $P(q, t)$ (cf. \cite{F2, P}).
\begin{prop} \label{prop1}
$(i)$ $P(q, t) : \D R \times \D R \to \D R$ is continuous.

$(ii)$ $P(q, t)$ is strictly decreasing in each variable separately.

$(iii)$ For fixed $q$ we have $\lim_{t \to +\infty} P(q, t)=-\infty$ and $\lim_{t \to -\infty} P(q, t)=+ \infty$.

$(iv)$ $P(q, t)$ is convex: if $q_1, q_2, t_1, t_2 \in \D R, \, a_1, a_2 \geq 0, a_1+a_2=1$, then
\[P(a_1q_1+a_2q_2, a_1t_1 +a_2 t_2) \leq a_1 P(q_1, t_1) +a_2 P(q_2, t_2).\]
\end{prop}

Now for fixed $q$, $P(q, t)$ is a continuous function of $t$. Its values range from $-\infty$ (when $t \to +\infty$) to $+\infty$ (when $t \to -\infty$). Therefore, by the intermediate value theorem there is a real number $\gb$ such that $P(q, \gb)=0$. The solution $\gb$ is unique, since $P(q, \cdot)$ is strictly decreasing. This defines $\gb$ implicitly as a function of $q$: for each $q$ there is a unique $\gb=\gb(q)$ such that $P(q, \gb(q))=0$.

The following proposition gives the well-known properties of the function $\gb(q)$ (cf. \cite{F2, P}).
\begin{prop}
Let $\gb=\gb(q)$ be defined by $P(q, \gb(q))=0$. Then

$(i)$ $\gb$ is a continuous function of the real variable $q$.

$(ii)$ $\gb$ is strictly decreasing: if $q_1<q_2$, then $\gb(q_1)> \gb(q_2)$.

$(iii)$ $\lim_{q \to -\infty}\gb(q)=+\infty$ and $\lim_{q \to +\infty}\gb(q)=-\infty$.

$(iv)$ $\gb$ is convex: if $q_1, q_2, a_1, a_2 \in \D R$ with $a_1, a_2 \geq 0$ and $a_1+a_2=1$, then
\[\gb(a_1q_1+a_2q_2) \leq a_1\gb(q_1)+a_2 \gb(q_2).\]
\end{prop}

The function $\gb(q)$ is sometimes denoted by $T(q)$ and called the \tit{temperature function}. A more general discussion of this function can be found in \cite{HJKPS}, where our $\gb(q)$ function corresponds to $-\gt(q)$ in their notation.

\begin{remark} If $q=0$, then $P(0, \gb(0))=0$, which implies
\[\lim_{k\to\infty} \frac 1 k \log \sum_{\gs\in \gO_k}\|\vp_\gs'\|^{\gb(0)}=0,\]
i.e., $\gb(0)$ gives the Hausdorff dimension $\te{dim}_\te{H}(E)$ of the cookie-cutter set $E$ (cf. \cite{F1}). Again
\[P(1, 0)=\lim_{k\to\infty} \frac 1 k \log \sum_{\gs\in \gO_k}\mu_h(J_\gs)=\lim_{k\to\infty} \frac 1 k \log 1=0,\]
and hence $\gb(1)=0$ (see Figure 1).
\end{remark}

\section{Main result}

The relationship between the quantization dimension function and the
temperature function $\gb(q)$ for the Gibbs-like measure $\mu_h$, where
the temperature function is the Legendre transform of the $f(\ga)$
curve (for the definitions of $f(\ga)$ and the Legendre transform see
\cite{F1}) is given by the following theorem which constitutes the
main result of the paper. For its graphical description see Figure 1.

\begin{theorem} \label{theorem}
Let $\mu_h$ be the Gibbs-like measure supported by the cookie-cutter set $E$. Then, for each $r
\in (0, +\infty)$ there exists a unique $\gk_r \in (0, +\infty)$ such that
\[\gk_r=\frac{\gb(q_r)}{1-q_r},\]
where we recall $\beta$ is the temperature function, i.e., $\gb(q_r)=rq_r$, and $\gk_r$ is the quantization dimension of order $r$ of the probability measure $\mu_h$. Moreover, the $\gk_r$-dimensional lower quantization coefficient is positive, i.e.,
$\liminf_n n V_{n,r}^{\frac{\gk_r}{r}}(\mu_h)>0$.
\end{theorem}

\begin{figure}
\includegraphics[width=4in]{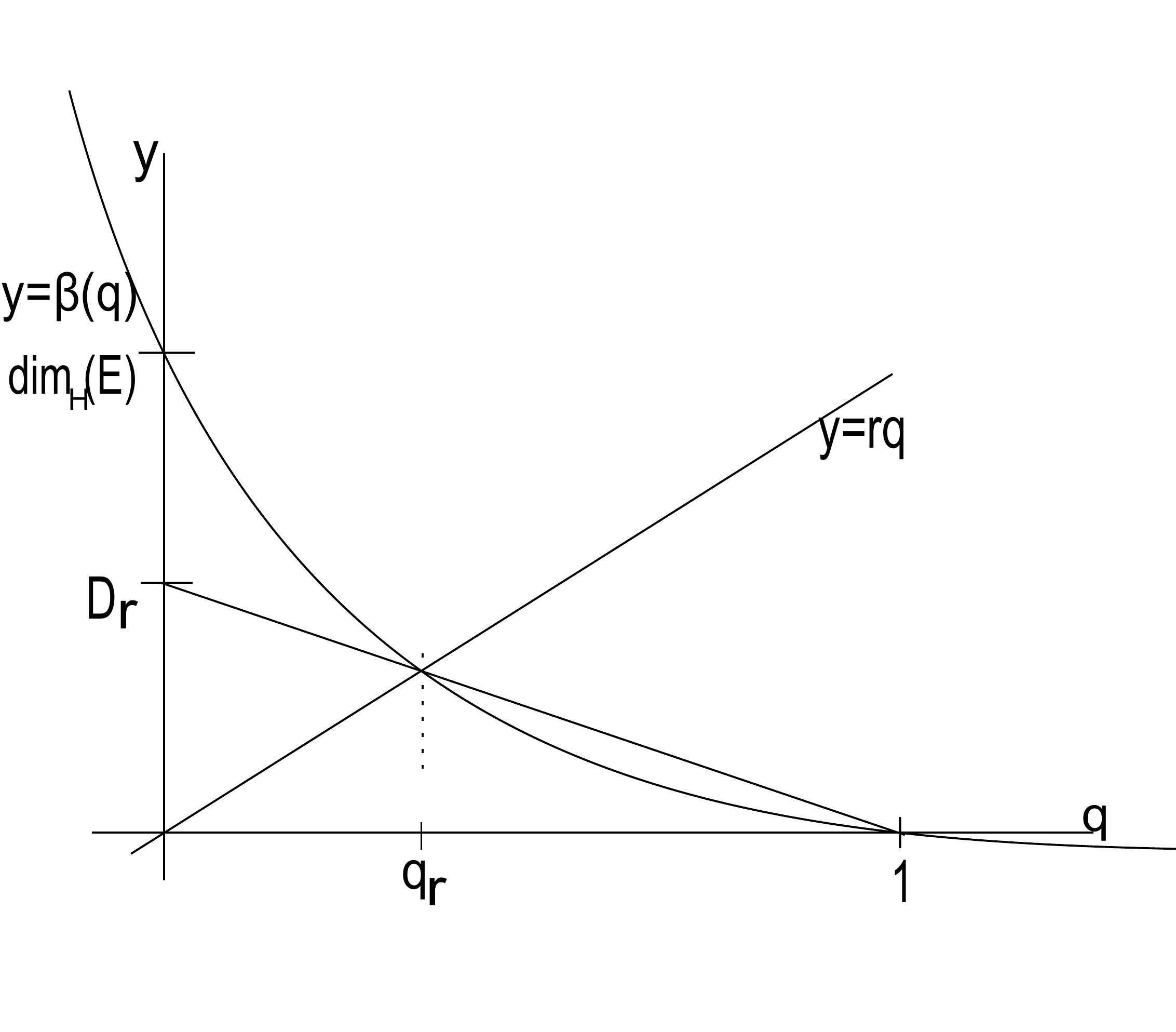}
\caption {\label{Fig1.} To determine $D_r$ first find the point of intersection of $y=\gb(q)$ and the line $y=rq$. Then $D_r$ is the $y$-intercept of the line through this point and the point (1, 0).}
\end{figure}

To prove the theorem we need several lemmas and propositions. Let us first state the following lemma, the proof is similar to \cite[Lemma~3.2]{R2}.
\begin{lemma}  \label{lemma4}
Let $0<r<+\infty$ be fixed. Then there exists exactly one number $\gk_r
\in (0, +\infty)$ such that \[ \lim_{k\to\infty} \frac 1 k \log \sum_{\gs\in \gO_k}\left(\mu_h(J_\gs) \|\vp_\gs'\|^r\right)^{\frac{\gk_r}{r+\gk_r}}=0.\]
\end{lemma}

Let us now prove the following lemma.

\begin{lemma} \label{lemma5}
Let $0<r<+\infty$ and $\gk_r$ be as in
Lemma~\ref{lemma4}. Then for $n\geq 1$, we have
\[(L \gx^{r} )^{-\frac{\gk_r}{r+\gk_r}}\leq \sum_{\gs\in \gO_n}\left(\mu_h(J_\gs)\|\vp_\gs'\|^r\right)^{\frac{\gk_r}{r+\gk_r}}\leq
(L\gx^{r} )^{\frac{\gk_r}{r+\gk_r}}.\]
\end{lemma}
\begin{proof} For $\gs\in \gO$ let us write $s_\gs=\mu_h(J_\gs)\|\vp_\gs'\|^r$. Then for
 $\gs \in \gO_n$ and $\gt\in \gO_p$ with $n, p\geq 1$, by Lemma~\ref{lemma345} and Lemma~\ref{lemma3}, we have $L^{-1} \gx^{-r}  s_\gs s_\gt \leq s_{\gs\gt}\leq L s_\gs s_\gt\leq  L \gx^r s_\gs s_\gt$. Since $r>0$, $L, \gx>1$, it is true that $ L^{-2} \gx^{-2r} s_\gs s_\gt\leq
 s_{\gs\gt}\leq L^{2}\gx^{2r} s_\gs s_\gt $. Then by the standard theory of sub-additive sequences, $\lim_{n\to \infty}n^{-1} \log\sum_{\gs\in \gO_n}s_{\gs}^t$
 exists for any $t\in \D R$. Then, proceeding as Lemma~3.3 in \cite {R2}, we obtain the lemma.

\end{proof}

We call $\gG \sci \gO$ a \tit{finite maximal antichain} if $\gG$ is a finite set of words in $\gO$, such that every sequence in $\gO_\infty$ is an extension of some word in $\gG$, but no word of $\gG$ is an extension of another word in $\gG$. By $|\gG|$ we denote the cardinality of $\gG$. Note that from the definition of $\gG$ it follows that finite maximal antichain does not contain the empty word $\es$ as all words are extension of $\es$.

Let us now state and prove the following lemma.
\begin{lemma}\label{lemma6}
Let $0<r<+\infty$ and  $\gk_r$ be as in Lemma~\ref{lemma4}. Then, for any finite maximal antichain $\gG$, we have

$(a)$  $L^{-1}\sum_{\gs\in \gG} \mu_h(J_\gs) \mu_h \circ \vp_{\gs}^{-1} \leq \mu_h \leq L \sum_{\gs\in \gG} \mu_h(J_\gs) \mu_h \circ \vp_{\gs}^{-1}$, and

$(b)$   $ (L\gx^{r} )^{-\frac{3\gk_r}{r+\gk_r}}\leq \sum_{\gs \in \gG}\left (\mu_h(J_\gs)\|\vp_\gs'\|^r\right)^{\frac{\gk_r}{r+\gk_r}}\leq  (L\gx^{r} )^{\frac{3\gk_r}{r+\gk_r}}.$

\end{lemma}

\begin{proof} $(a)$ Let $M=\max\set{|\gs| : \gs \in \gG}$. Note that the Borel $\gs$-algebra on $E$ is generated by the set $\set{J_\gs : \gs \in \gO}$ of all
basic intervals. For any two basic intervals either they are disjoint or one is contained in the other. Hence, it is enough to prove that for any $\go\in \gO_n$ with
$n\geq M$,
\[L^{-1}\sum_{\gs\in \gG} \mu_h(J_\gs) \mu_h \circ \vp_{\gs}^{-1}(J_\go) \leq \mu_h(J_\go) \leq L \sum_{\gs\in \gG} \mu_h(J_\gs) \mu_h \circ \vp_{\gs}^{-1}(J_\go),\]
which follows in the similar lines as the proof of Lemma~3.4 $(a)$ in \cite {R2}.

$(b)$ As $\gG$ is a finite maximal antichain, there exists a finite sequence of positive integers $n_1<n_2<\cdots <n_K$ such that
\[\gG=\gG_1\uu\gG_2\uu\cdots \uu \gG_K,\]
where $\gG_j=\set{\gs \in \gG : |\gs|=n_j}$ for all $1\leq j\leq K$. Let $M$ be a positive integer such that $M\geq n_K$. Then by Lemma~\ref{lemma5}, we have
\begin{align*}
&\sum_{\gs \in \gG}\left(\mu_h(J_\gs)\|\vp_\gs'\|^r\right)^{\frac{\gk_r}{r+\gk_r}}\\
&\geq \sum_{j=1}^K \sum_{\gt\in\gG_j}\left(\mu_h(J_\gt)\|\vp_\gt'\|^r\right)^{\frac{\gk_r}{r+\gk_r}}(L\gx^{r} )^{-\frac{\gk_r}{r+\gk_r}}\sum_{\gs\in \gO_{M-n_j}}\left(\mu_h(J_\gs)\|\vp_\gs'\|^r\right)^{\frac{\gk_r}{r+\gk_r}} \\
& \geq (L\gx^{r} )^{-\frac{\gk_r}{r+\gk_r}}\sum_{j=1}^K \sum_{\gt\in\gG_j}\left(\mu_h(J_\gt)\|\vp_\gt'\|^r\right)^{\frac{\gk_r}{r+\gk_r}}\sum_{\begin{subarray}{c} {\gs \in \gO_{M-n_j}}\\{\gt\prec\gs} \end{subarray}}\left(\mu_h(J_\gs)\|\vp_\gs'\|^r\right)^{\frac{\gk_r}{r+\gk_r}}\\
& = (L\gx^{r} )^{-\frac{\gk_r}{r+\gk_r}}\sum_{j=1}^K \sum_{\gt\in\gG_j}\sum_{\begin{subarray}{c} {\gs \in \gO_{M-n_j}}\\{\gt\prec\gs} \end{subarray}}\left(\mu_h(J_\gt)\|\vp_\gt'\|^r\right)^{\frac{\gk_r}{r+\gk_r}} \left(\mu_h(J_\gs)\|\vp_\gs'\|^r\right)^{\frac{\gk_r}{r+\gk_r}}\\
& \geq (L\gx^{r} )^{-\frac{2\gk_r}{r+\gk_r}}\sum_{j=1}^K \sum_{\gt\in\gG_j}\sum_{\begin{subarray}{c} {\gs \in \gO_{M-n_j}} \end{subarray}} \left(\mu_h(J_{\gt\gs})\|\vp_{\gt\gs}'\|^r\right)^{\frac{\gk_r}{r+\gk_r}}\\
&= (L\gx^{r} )^{-\frac{2\gk_r}{r+\gk_r}}\sum_{\gs\in \gO_M}\left(\mu_h(J_\gs)\|\vp_\gs'\|^r\right)^{\frac{\gk_r}{r+\gk_r}}\\
&\geq (L\gx^{r} )^{-\frac{3\gk_r}{r+\gk_r}}.
\end{align*}
Similarly, we have $\sum_{\gs \in \gG}\left(\mu_h(J_\gs)\|\vp_\gs'\|^r\right)^{\frac{\gk_r}{r+\gk_r}}\leq  (L\gx^{r} )^{\frac{3\gk_r}{r+\gk_r}}$. Hence the lemma.
\end{proof}
Let us now give the following lemma.
\begin{lemma}  Let $x, y \in J$ and $\gs \in \gO$. Then
\[\gx^{-1} \|\vp_\gs'\| d(x, y) \leq d(\vp_\gs(x), \vp_\gs(y)) \leq \|\vp_\gs'\| d(x, y),\]
where $\gx$ is the constant of Lemma~\ref{lemma1111}.
\end{lemma}
\begin{proof}
By mean value theorem, for any $x, y\in J$ there exists some $w\in (x, y)$ such that
\[d(\vp_\gs(x), \vp_\gs(y))=|\vp_\gs'(w)| d(x, y), \]
and so, by Proposition~\ref{prop32}, for any $z \in J$,
\[\gx^{-1} |\vp_\gs'(z)|d(x, y) \leq d(\vp_\gs(x), \vp_\gs(y))=|\vp_\gs'(w)| d(x, y)\leq \|\vp_\gs'\|d(x, y).\]
Now take the supremum over all $z\in J$, and then
\[\gx^{-1} \|\vp_\gs'\| d(x, y) \leq d(\vp_\gs(x), \vp_\gs(y)) \leq \|\vp_\gs'\|d(x, y),\]
to obtain the assertion of the lemma.
\end{proof}

Using the above lemma and Lemma~\ref{lemma6} $(a)$, and the parallel lines as Lemma~3.5 in \cite{R2} the following lemma can be proved.
\begin{lemma} \label{lemma7}
Let $\gG\sci \gO $ be a finite maximal antichain, $n \in \D N$ with $n \geq |\gG|,$ and $0<r <+\infty$. Then
$V_{n, r}(\mu_h)\leq \inf \left\{L \sum_{\gs\in \gG}\mu_h(J_\gs) \|\vp_\gs'\|^rV_{n_\gs, r}( \mu_h) : 1 \leq n_{\gs}, \ \sum_{\gs \in \gG}  n_{\gs} \leq n \right\}.$\end{lemma}

Using the above lemma and similar lines as Proposition~3.6 in \cite{R2}, the following proposition can be proved.

\begin{prop}\label{prop5}
\tit{Let $0<r<+\infty$ and $\gk_r$ be as in Lemma~\ref{lemma4}. Then $\limsup_n n V_{n, r}^{\frac{\gk_r}{r}}(\mu_h) <+ \infty.$}
\end{prop}

\begin{note}

We say that the cookie-cutter mapping $f$ satisfies the \tit{open set condition} (OSC) if the corresponding set
of branch inverses $\set{\vp_1, \vp_2, \cdots, \vp_N}$ satisfies the open set condition: there exists a bounded nonempty
open set $U \sci J$ (in the topology of $J$) such that $\vp_{j}(U) \sci U$, and $\vp_i(U) \ii \vp_j(U)=\es$ for $1\leq i\neq j\leq N$. Furthermore,
$f$ satisfies the \tit{strong open set condition} (SOSC) if $U$ can be chosen such that $U\ii E\neq \es$. Note that we can choose $U=(0, 1)$ and so,
because of the net properties of the basic intervals of the cookie-cutter set, it follows that the cookie-cutter mapping $f$  satisfies the strong open set condition.
\end{note}

Parallel to Lemma~3.7 in \cite{R2} the following lemma can be proved.
\begin{lemma} \label{lemma105}
Let $\gG \sci \gO$ be a finite maximal antichain. Then there exists $n_0=n_0(\gG)$ such that for every $n\geq n_0$ there exists a set of positive integers $\set {n_\gs: =n_\gs(n)}_{\gs \in \gG}$ such that $\sum_{\gs \in \gG}n_\gs\leq n$ and
\[u_{n, r}(\mu_h) \geq \left(L\gx^r\right)^{-1} \sum_{\gs \in \gG}\mu_h(J_\gs) \|\vp_\gs'\|^r  u_{n_\gs, r}(\mu_h).\]
\end{lemma}

Let us now prove the following proposition, which shows that the $\gk_r$-dimensional lower quantization coefficient of order $r$ of the probability measure $\mu_h$ is positive.
\begin{prop} \label{prop6}
\tit{Let $\mu_h$ be the Gibbs-like measure, and let $0<r <+\infty$. Moreover, let $\gk_r$ be as in Lemma~\ref{lemma4}. Then $\liminf_n n
V_{n, r}^{\frac{\gk_r}{r}}(\mu_h)>0.$}
\end{prop}
\begin{proof}
Let $\gG$ be a finite maximal antichain. By Lemma~\ref{lemma105}, we have $n_0$ and for $n\geq n_0$ the numbers $\set{n_\gs:=n_\gs(n)}_{\gs \in \gG}$ which satisfy the conclusion of the lemma. Set
$c=\min \set{n^{r/\gk_r} u_{n, r}(\mu_h) : n \leq n_0}$. Clearly each $u_{n, r}(\mu_h)>0$ and hence $c>0$. Suppose $n \geq n_0$ and $k^{r/\gk_r} u_{k, r}(\mu_h) \geq c$ for all $k<n$. Hence, using Lemma~\ref{lemma105}, we have
\begin{align*}
n^{r/\gk_r} u_{n, r}(\mu_h) & \geq  n^{r/\gk_r} (L\gx^r)^{-1}\sum_{\gs \in \gG} \mu_h(J_\gs)\|\vp_\gs'\|^r u_{n_\gs, r}(\mu_h)\\
&=n^{r/\gk_r} (L\gx^r)^{-1}  \sum_{\gs \in \gG} \mu_h(J_\gs)\|\vp_\gs'\|^r(n_\gs(n))^{-r/\gk_r}(n_\gs(n))^{r/\gk_r} u_{n_\gs, r}(\mu_h)\\
& \geq c (L\gx^r)^{-1} \sum_{\gs \in \gG} \mu_h(J_\gs)\|\vp_\gs'\|^r \left(\frac{n_\gs(n)}{n}\right)^{-r/\gk_r}.
\end{align*}
Using H\"older's inequality (with exponents less than 1), we have
\[n^{r/\gk_r} u_{n, r}(\mu_h)\geq c (L\gx^r)^{-1} \left( \sum_{\gs \in \gG}\left(\mu_h(J_\gs)\|\vp_\gs'\|^r\right)^{\gk_r/(r+\gk_r)}\right)^{(1+r/\gk_r)}\left(\sum_{\gs \in \gG}(\frac{n_\gs(n)}{n})^{(-r/\gk_r)(-\gk_r/r)}\right)^{-r/\gk_r}.\]
By Lemma~\ref{lemma6} $(b)$, and the fact that $\sum_{\gs \in\gG} n_\gs(n)
\leq n$, we have
\[n^{r/{\gk_r}} u_{n, r}(\mu_h)\geq  c (L\gx^r)^{-1}(L\gx^r)^{-3}.\]
Therefore, by induction,
\[\liminf_n n u_{n, r}^{\gk_r/r}(\mu_h) \geq [c(L\gx^r)^{-4}]^{\gk_r/r}>0, \te{ i.e., } \liminf_n nV_{n, r}^{\frac {\gk_r}{r}}(\mu_h)>0,\]
and thus the proposition is yielded.
\end{proof}

\subsection*{Proof of Theorem~\ref{theorem}}

Note that $e_{n, r}=V_{n, r}^{\frac 1 r}$, and by Proposition 11.3 of \cite{GL1}, we know:

$(a)$ If $0\leq t<\ul D_r<s$ then \[\lim_{n\to \infty} ne^t_{n, r}=+\infty \te{ and } \liminf_{n\to \infty} ne^s_{n, r}=0.\]

$(b)$ If $0\leq t<\ol D_r<s$ then
\[\limsup_{n\to \infty} ne^t_{n, r}=+\infty \te{ and } \lim_{n\to \infty} ne^s_{n, r}=0.\]

From $(a)$ and Proposition~\ref {prop6}, we have  $\gk_r\leq \ul D_r$. From $(b)$ and Proposition~\ref{prop5},
we have $\ol D_r \leq \gk_r$. Hence, $\gk_r\leq \ul D_r\leq \ol D_r \leq \gk_r$, i.e., the quantization dimension $D_r$ exists and $D_r=\gk _r.$ Note that if $q_r=\frac{\gk_r}{r+\gk_r}$,  by Lemma~\ref{lemma4}, we have $\gb(q_r)=rq_r$, and then $D_r=\frac{\gb(q_r)}{1-q_r}$. Moreover, by Proposition~\ref{prop6}, we have $\liminf_n n V_{n, r}^{\gk_r/r}(\mu_h)>0$.  Thus the proof of the theorem is complete.\\

\end{document}